\documentclass[10pt; a4paper]{amsart} 
\usepackage{amscd, amsfonts, amsmath, amssymb, amsthm, fixmath, mathrsfs}
\usepackage[all]{xy} 
\makeindex

\theoremstyle{definition} 
\newtheorem{Unity}{Unity}[section] 
\newtheorem*{Definition*}{Definition} 
\newtheorem{Definition}[Unity]{Definition} 

\theoremstyle{plain} 
\newtheorem*{Theorem*}{Theorem}
\newtheorem{Theorem}[Unity]{Theorem}
\newtheorem{Proposition}[Unity]{Proposition}

\newtheorem{Lemma}[Unity]{Lemma}

\theoremstyle{remark} 
\newtheorem*{Remark*}{Remark}

\numberwithin{Unity}{section}

\newcommand{\p}{\mathfrak{p}}
\newcommand{\q}{\mathfrak{q}}
\newcommand{\m}{\mathfrak{m}}
\newcommand{\V}{\mathrm{V}}

\newcommand{\hgt}{\mathrm{ht}}

\newcommand{\Max}{\mathrm{Max}}
\newcommand{\Ann}{\mathrm{Ann}}
\newcommand{\Att}{\mathrm{Att}}
\newcommand{\Cos}{\mathrm{Cos}}
\newcommand{\Hom}{\mathrm{Hom}}

\newcommand{\Spec}{\mathrm{Spec}}
\newcommand{\Supp}{\mathrm{Supp}}
\newcommand{\Ndim}{\mathrm{Ndim}}

\newcommand{\Cdim}{\mathrm{Cdim}}

\newcommand{\cograde}{\mathrm{cograde}}

\begin{document}
\title{Dual Bass Numbers and Co-Cohen Macaulay Modules}
\author{Lingguang Li}
\address{Department of Mathematics, Tongji University, Shanghai, P. R. China\\ School of Mathematical Sciences, Fudan University, Shanghai, P. R. China}
\email{LG.Lee@amss.ac.cn}
\begin{abstract} In this paper, we give a characterization of co-Cohen Macaulay modules by vanishing properties of the dual Bass numbers of modules. In addition, we show that the co-localization of co-Cohen Macaulay modules preserves co-Cohen Macaulayness under a certain condition.
\end{abstract}
\maketitle

\section{Introduction}

Let $R$ be a Noetherian ring, $M$ an $R$-module, $\p\in\Spec~R$. H. Bass defined so called Bass numbers $\mu_{i}(\p,M)$ by using the minimal injective resolution of $M$ for all integers $i\geq 0$, and proved that $$\mu_{i}(\p,M)=\dim_{k(\p)}\textmd{Ext}_{R_{\p}}^{i}(k(\p),M_{\p}).$$
If M is a finitely generated $R$-module, the Betti numbers $\beta_{i}(\p,M)$ is defined by using the minimal free resolution of $M_{\p}$ for all $i\geq
0$, and we have $$\beta_{i}(\p,M)=\dim_{k(\p)}\textmd{Tor}^{R_{\p}}_{i}(k(\p),M_{\p}).$$ In addition, E. Enochs and J. Z. Xu defined the dual Bass numbers $\pi_{i}(\p,M)$ by using the minimal flat resolution of $M$ for all $i\geq 0$ and showed that $$\pi_{i}(\p,M)=\dim_{k(\p)}\textmd{Tor}^{R_{\p}}_{i}(k(\p),\Hom_{R}(R_{\p},M))$$ for any cotorsion $R$-module $M$ in \cite{EnochsXu97}. J. Z. Xu characterized Gorenstein rings and strongly cotorsion modules by vanishing properties of $\pi_{i}(\p,M)$. The author have studied the vanishing properties of dual Bass numbers in \cite{Li10}.

A finitely generated $R$-module $M$ over a Noetherian ring $R$ is called \emph{Cohen Macaulay} if $\dim_{R_{\p}}M_{\p}=\mathrm{grade}_{R_{\p}}M_{\p}$ for any $\p\in\Supp_RM$. Cohen Macaulay modules over Noetherian rings are important objects in commutative algebra and algebraic geometry. In duality, Z. M. Tang and H. Zakeri introduced the concept of co-Cohen Macaulay modules and studied the properties of this in \cite{TangZakeri94} and \cite{Tang96}. An elementary and important property of Cohen Macaulay modules is that the localization preserves the Cohen Macaulayness. The dual question for Artinian modules is to ask whether the co-localization of co-Cohen Macaulay modules preserves the co-Cohen Macaulayness. We show that this statement is true under a certain conditions (Proposition \ref{Prop:Cograde-Cdim}).

In addition, we give a characterization of co-Cohen Macaulay modules by vanishing properties of dual Bass numbers which have relation to the maximal length of co-regular sequence (Theorem \ref{Thm:Co-CM}). This is dual to the theory of Cohen Macaulay modules over Noetherian rings.

\section{Preliminaries}

The concept of Krull dimension (Kdim) for Artinian modules was introduced by R. N. Roberts in \cite{Roberts75}. Later, Kirby \cite{Kirby90} changed the terminology of Roberts and referred it to Noetherian dimension (Ndim) to avoid any confusion with well-known Krull dimension defined for finitely generated modules. Let $R$ be a ring, $M$ a $R$-module. The Noetherian dimension of $M$, denoted by $\Ndim_RM$, is defined inductively as follows: when $M=0$, put $\Ndim_RM=-1$. Then by induction, for an integer $d\geq 0$, we put $\Ndim_RM=d$, if $\Ndim_RM=d$ is false and for every ascending chain $M_0\subseteq M_1\subseteq\cdots$ of submodules of $M$, there exists a positive integer $n_0$ such that $\Ndim_R(M_{n+1}/M_n)<d$ for all $n>n_0$. Therefore $\Ndim_RM=0$ if and only if $M$ is a non-zero Noetherian module.

L. Melkersson and P. Schenzel introduced the co-localization of modules in \cite{MelSch95}. Let $R$ be a ring, $S\subseteq R$ a multiplicative set, and $M$ an $R$-module. The $R_S$-module $\Hom_{R}(R_S,M)$ is called the \emph{co-localization} of $M$ with respect to $S$, and defined the \emph{co-support} of $M$ by $\Cos_RM=\{~\p\in\Spec~R~|~\Hom_{R}(R_{\p},M)\neq 0~\}$. In \cite{Li10}, we define the \emph{co-dimension} of $M$ as $$\Cdim_RM=\sup\{~\dim R/\p~|~\p\in\Cos_RM~\}.$$

If $R$ is a Noetherian ring and $M$ an Artinian $R$-module, then $\Att_RM$, $\Cos_RM$, $\Ann_RM$ have the same minimal elements by \cite[Corollary 4.3]{CuongNhan02i}.

Let $R$ be a ring (not necessarily Noetherian), $S\subseteq R$ a multiplicative set, $M$ an Artinian $R$-module. Then we have
\begin{eqnarray*}
\Cdim_{S^{-1}R}\Hom_{R}(S^{-1}R,M)&=&\sup\{\dim S^{-1}R/S^{-1}\p|~\p\in\Cos_RM,\p\cap S=\emptyset\}\\
&=&\sup\{\dim S^{-1}R/S^{-1}\p~|~\p\in\Att_RM,~\p\cap S=\emptyset\}.
\end{eqnarray*}
Let $\p\in\Cos_RM$, we denote$$\hgt_M\p =\sup\{~n~|~\p_{0}\subsetneq\p_{1}\subsetneq\cdots\subsetneq \p_{n},~\p_{i}\in\Cos_RM~\text{for}~i=0,1,\cdots,n~\}.$$ It is obvious that $\hgt_M\p=\Cdim_{R_{\p}}\Hom_{R}(R_{\p},M)$.

We expect that $\Ndim_RM=\Cdim_RM$ for any Artinian module $M$. Unfortunately, this equality may not hold in general. In fact, there exists an Artinian module $M$ over a Noetherian local ring $(R,\m)$ such that $\Ndim_RM<\Cdim_RM$, See \cite[Example 4.1]{CuongNhanN02ii}.
However, we have the following proposition.

N. T. Cuong and N. T. Dung showed that the above condition does not hold for any Artinian modules, and they also gave some sufficient conditions for that in \cite{CDN}.

\cite[Proposition 2.1]{CDN} Let $(R,\m)$ be a Noetherian local ring, $M$ an Artinian $R$-module, if one of the following cases
happens:
\begin{itemize}
\item[$(i)$.] $R$ is complete with respect to $m$-adic topology.
\item[$(ii)$.] $M$ contain a submodule which is isomorphic to the injective hull of $R/\m$.
\end{itemize}
Then $\Ann_R(0:_M\p)=\p$ for any $\p\in\V(\Ann_RM)$.

\begin{Definition} A Noetherian ring $R$ is called a \emph{$U$ ring}, if for any Artinian $R$-module $M$ such that $\Ann_R(0:_M\p)=\p$ for any $\p\in\V(\Ann_RM)$.
\end{Definition}

Hence, any complete Noetherian local rings are $U$ rings.

Let $R$ be a ring, $M$ an Artinian ring. R. Y. Shorp in \cite{Sharp89} showed that $\Supp_RM$ is a finite subset of $\Max~R$, and if $\Supp_RM=\{\m_1,\cdots,\m_s\}$, then $M=M_{1}\oplus\cdots\oplus M_{s}$, where $M_i=\{x\in M|m_i=\sqrt{\Ann_Rx}\}$.

\section{Co-localization of Co-Cohen Macaulay Modules}

Cohen-Macaulay modules over Noetherian rings are important objects in commutative algebra. In duality, Z. Tang and H. Zakeri introduced the co-Cohen Macaulay modules over local rings in \cite{TangZakeri94} and it was generalized to general rings by Z. Tang in \cite{Tang96}.

\begin{Definition}\cite[Definition 5.3]{Tang96}
Let $R$ be a ring. An Artinian $R$-module $M$ is called a \emph{co-Cohen-Macaulay $R$-module}, if $\cograde_{R}(J(M),M)=\Ndim_RM$, where $\cograde_{R}(J(M),M)$ is the common length of any maximal $M$ co-regular sequence contained in $J(M)$.
\end{Definition}

\begin{Proposition}\label{Prop:CoCMDecomp}
Let $R$ be a ring, $M$ an Artinian $R$-module. If $M=M_1\oplus\cdots\oplus M_s$ and $\Supp_RM=\{\m_1,\cdots,\m_s\}$, where $M_i=\{x\in M|m_i=\sqrt{\Ann_Rx}\}$. Then $M$ is a co-Cohen-Macaulay $R$-module if and only if $M_{i}$ is a co-Cohen-Macaulay $R$-module and $\Ndim_RM=\Ndim_RM_{i}$ for $i=1,2,\cdots,s$.
\end{Proposition}

\begin{proof}
Notice that $\Ndim_RM=\textmd{max}\{\Ndim_RM_{i}|i=1,2,\cdots,s\}$ by \cite[Proposition 1]{Roberts75}.
Since $M$ is a co-Cohen Macaulay $R$-module, and $J(M_{i})=\m_{i}$ for any $1\leq i\leq s$, by \cite[Lemma 5.1 and Proposition 5.2]{Tang96}, we get
$$\Ndim_RM=\Ndim_RM_{i}=\cograde_{R}(\m_{i},M_{i})=\cograde_{R}(J(M),M)$$ for $i=1,\cdots,s$. Hence, $M_{i}$ is a co-Cohen Macaulay $R$-module and $\Ndim_RM=\Ndim_RM_{i}$ for any $i=1,2,\cdots,s$. It is similar to prove the other side.
\end{proof}

It is well-known that operation of localization preserves the Cohen Macaulayness of Cohen Macaulay modules over Noetherian rings. Since co-localization is dual to localization, hence there is a natural question: Whether co-localization preserves the co-Cohen Macaulayness of co-Cohen Macaulay modules? The main content of this section is to study the properties of co-localization of co-Cohen Macaulay modules, and give a partial affrmative answer to above question.

\begin{Proposition}\label{Prop:CdCoLoc}
Let $R$ be a Noetherian ring, $M$ an Artinian $R$-module such that $\Ann_R(0:_M\p)=\p$ for any $\p\in\V(\Ann_RM)$. If $\p\in\Cos_RM$,and $\frac{x}{s}\in \p R_{\p}$ is a $\Hom_{R}(R_{\p},M)$ co-regular element. Then
$$\Cdim_{R_{\p}}(0:_{\tiny{\Hom_{R}(R_{\p},M)}}\frac{x}{s})=\Cdim_{R_{\p}}\Hom_{R}(R_{\p},M)-1.$$
\end{Proposition}

\begin{proof}Since $(0:_{\tiny{\Hom_{R}(R_{\p},M)}}\frac{x}{s})\cong\Hom_{R}(R_{\p},0:_Mx)$, we have
$$\Cdim_{R_{\p}}(0:_{\tiny{\Hom_{R}(R_{\p},M)}}\frac{x}{s})=\hgt(\p/\Ann_R(0:_Mx)).$$
Moreover, by \cite[Lemma 4.3]{Li10} we have $\V(\Ann_R(0:_Mx))=\V(\Ann_RM,x).$
Hence $\sqrt{\Ann_R(0:_Mx)}=\sqrt{(\Ann_RM,x)}$. Since $\Cdim_{R_{\p}}\Hom_{R}(R_{\p},M)=\hgt(\p/\Ann_RM)$,
we only need to show that $\hgt(\p/(\Ann_RM,x))=\hgt(\p/\Ann_RM)-1$.

Since $\frac{x}{s}\in\p R_{p}$ is a $\Hom_{R}(R_{\p},M)$ co-regular element, we get $$\frac{x}{s}\in\p R_{\p}-\bigcup\limits_{\tiny{\q\subseteq \p\atop \q\in\Att_RM}}\q R_{\p}.$$ It follows that $x\in \p-\bigcup\limits_{\tiny{\q\subseteq\p\atop\q\in\Att_RM}}\q$. Thus
$\hgt(\p/(\Ann_RM,x))\leq\hgt(\p/\Ann_RM)-1$. On the other hand, we get $\hgt(\p/(\Ann_RM,x))\geq\hgt(\p/\Ann_RM)-1$
by the theory of system of parameters. Hence, the result follows.
\end{proof}

The following Proposition shows that co-localization preserves co-Cohen Macaulayness under a certain conditions.

\begin{Proposition}\label{Prop:Cograde-Cdim}
Let $R$ be a Noetherian ring, $M$ a co-Cohen Macaulay $R$-module, such that $\Ann_R(0:_M\p)=\p$ for any $\p\in\V(\Ann_RM)$. Then for any $\p\in\Cos_RM$, $$\cograde_{R}(\p,M)=\cograde_{R_{\p}}\Hom_{R}(R_{\p},M)=\Cdim_{R_{\p}}\Hom_{R}(R_{\p},M).$$
\end{Proposition}

\begin{proof}Let $\p\in\Cos_RM$, then we have $\cograde_{R_{\p}}\Hom_{R}(R_{\p},M)<\infty$, and
$$\cograde_{R}(\p,M)\leq\cograde_{R_{\p}}\Hom_{R}(R_{\p},M) \leq\Cdim_{R_{\p}}\Hom_{R}(R_{\p},M).$$
We only need to show that $\cograde_{R}(\p,M)=\Cdim_{R_{\p}}\Hom_{R}(R_{\p},M)$.

Let $n=\cograde_{R}(\p,M)$, we use induction on $n$. For the case $n=0$, there exists ${\frak Q}\in\Att_RM$ such that $\p\subseteq{\frak Q}$. Since $\p\in\Cos_RM$, then for any $\q\in\Att_RM$ with $\q\subseteq\p$, we have $\q=\p={\frak Q}$ by \cite[Proposition 5.2]{Tang96} and co-Cohen Macaulayness of $M$. Thus $\Att_{R_{\p}}\Hom_{R}(R_{\p},M)=\{\p R_{\p}\}$. Hence, $\Cdim_{R_{\p}}\Hom_{R}(R_{\p},M)=0$.

Suppose that $n>0$ and the Proposition holds for $n-1$. Let $x\in\p$ be a $M$ co-regular element, then $\cograde_{R}(\p,0:_Mx)=n-1$. By induction hypothesis, we have $\cograde_{R}(\p,0:_Mx)=\Cdim_{R_{\p}}\Hom_{R}(R_{\p},0:_Mx)$. On the other hand, by \cite[Proposition 3.4]{Li10}, $\frac{x}{1}\in \p R_{\p}$ is a $\Hom_{R}(R_{\p},M)$-quasi co-regular element. Since $(0:_{\tiny{\Hom_{R}(R_{\p},M)}}\p R_{\p})\neq 0$, we know that $\frac{x}{1}$ is a $\Hom_{R}(R_{\p},M)$ co-regular element. Then we have $\Cdim_{R_{\p}}\Hom_{R}(R_{\p},0:_Mx)=\Cdim_{R_{\p}}\Hom_{R}(R_{\p},M)-1$ by Proposition \ref{Prop:CdCoLoc}. Hence $\cograde_{R}(\p,M)=\Cdim_{R_{\p}}\Hom_{R}(R_{\p},M)$.
\end{proof}

To obtain the main result of this section, we first prove the following Lemma.

\begin{Lemma}\label{CosSupp}
Let $R$ be a Noetherian ring, $M$ an Artinian $R$-module. Then $$\Cos_RM\cap\V(J(M))=\Supp_RM.$$
\end{Lemma}

\begin{proof}
Let $M=N_1+N_2+\cdots+N_n$ be a minimal secondary presentation of $M$ such that $\Att_RM=\{\p_1,\cdots,\p_n\}$, where $\p_i=\sqrt{\Ann_RN_i}$, for $i=1,2,\cdots,n$. On the other hand, we have decomposition $M=M_{1}\oplus\cdots\oplus M_{s}$, $\Supp_RM=\{\m_1,\cdots,\m_s\}$, where $M_i=\{x\in M|m_i=\sqrt{\Ann_Rx}\}$. Assume that there exists $\m_j\in\Supp_RM-\Cos_RM$, then $\p_i\nsubseteq\m_j$ for any $1\leq i\leq n$. It follows that $\bigcap_{i=1}^n\p_i\nsubseteq\m_j$. Let $x\in\bigcap_{i=1}^n\p_i-\m_j$, then $x^l\cdot M=0$ for some integer $l$. Let $u$ be a non-zero element in $M_i$, then there exists $n\geq 0$ such that $\m_{i}^{n}u=0$. Thus $(\m_{i}^{n},x^l)u=0$. Since $\m_{i}$ is a maximal ideal, we get $(\m_{i}^{n},x^l)=R$. Thus $u=0$. This induces a contradiction. Hence, $\Supp_RM\subseteq\Cos_RM$. This implies $\Cos_RM\cap\V(J(M))=\Supp_RM$.
\end{proof}

\begin{Proposition}\label{Prop:}
Let $R$ be a Noetherian ring, $M$ an Artinian $R$-module such that $\Ann_R(0:_M\p)=\p$ for any $\p\in\V(\Ann_RM)$. Then $$\max\{\cograde_{R}(\p,M)|\p\in\Cos_RM\}=\max\{\cograde_{R}(\m,M)|\m\in\Supp_RM\}$$
\end{Proposition}

\begin{proof}
Since $\Ann_R(0:_M\p)=\p$ for any $\p\in\V(\Ann_RM)$, we have $\cograde_{R}(\p,M)<\infty~\text{for any}~\p\in\Cos_RM$.
By Lemma \ref{CosSupp}, we only need to show that for any $\p\in\Cos_RM$, there exists $\m\in\Supp_RM$ such that $\cograde_{R}(\p,M)\leq\cograde_{R}(\m,M)$.

If $x$ is a $M$ co-regular element, then we have $x\in\bigcup\limits_{\tiny{\m_{i}\in\Supp_RM}}\m_{i}$ by proof of Proposition \ref{Prop:Coreg-Ndim}. Let $\p\in\Cos_RM$, and $x_1,x_2,\cdots,x_{t}$ be a maximal $M$ co-regular sequence contained in $\p$, then there exists $\m\in\Supp_RM$ such that $x_i\in \m$ for $i=1,2,\cdots,t$. Otherwise, for any $1\leq k\leq s$, there exists $x_{i_{k}}\not\in\m_{k}$ for some $1\leq i_{k}\leq t$. Then $$(0:_M(x_1,x_2,\cdots,x_{t}))=(0:_{M_{1}}(x_1,x_2,\cdots,x_{t}))\oplus\cdots\oplus(0:_{M_{s}}(x_1,x_2,\cdots,x_{t}))=0.$$
This is a contradiction. Hence $\cograde_{R}(\p,M)\leq\cograde_{R}(\m,M)$ for some $\m\in\Supp_RM$.
\end{proof}

The following Theorem is the main result of this section. We give a characterization of co-Cohen Macaulay modules by vanishing
properties of dual Bass numbers.

\begin{Theorem}\label{Thm:Co-CM}
Let $R$ be a Noetherian ring, $M$ an Artinian $R$-module such that $\Ann_R(0:_M\p)=\p$
for any $\p\in\V(\Ann_RM)$. Let $\Max~\Cos_RM$ be the maximal elements of $\Cos_RM$. Then the following conditions are equivalent:
\begin{itemize}
\item[$(1)$.] $M$ is a co-Cohen Macaulay $R$-module.
\item[$(2)$.] $\cograde_{R_{\p}}\Hom_{R}(R_{\p},M)=\hgt_{\tiny{M}}\p$ for any $\p\in\Cos_RM$, and $\hgt_M\m=\Cdim_RM$ for any $\m\in\Supp_RM$.
\item[$(3)$.] $\pi_{i}(\p,M)=0$ for any $0\leq i<\hgt_M\p$ and for any ${\p}\in\Cos_RM$, and $\hgt_M\m=\Cdim_RM$ for any $\m\in\Supp_RM$.
\item[$(4)$.] $\cograde_{R}(\m,M)=\hgt_M\m=\Cdim_RM$ for any $\m\in\Supp_RM$.
\item[$(5)$.] $\pi_{i}(\m,M)=0$ for any $0\leq i<\hgt_M\m$, and $\hgt_M\m=\Cdim_RM$ for any $\m\in\Supp_RM$.
\end{itemize}
\end{Theorem}

\begin{proof} Notice that for any $\p\in\Cos_RM$, $\cograde_{R_{\p}}\Hom_{R}(R_{\p},M)<\infty$ by \cite[Lemma 4.3]{Li10}. Moreover,
$\Cdim_{R_{\p}}\Hom_{R}(R_{\p},M)=\hgt_M\p$, and $$\cograde_{R}(\m,M)=\cograde_{R_{\m}}\Hom_{R}(R_{\m},M)$$ for any $\m\in\Supp_RM$. Hence, by \cite[Proposition 5.2 and Proposition 5.9]{Li10}, we have $(2)\Leftrightarrow(3)$ and $(4)\Leftrightarrow(5)$. On the other hand, $(2)\Rightarrow(4)\Rightarrow(1)$ is obvious.

$(1)\Rightarrow(2)$. From Proposition \ref{Prop:Cograde-Cdim} we get $\cograde_{R_{\p}}\Hom_{R}(R_{\p},M)=\hgt_{\tiny{M}}\p$ for any $\p\in\Cos_RM$. Since $\cograde_{R}(\m,M)=\cograde_{R_{\m}}\Hom_{R}(R_{\m},M)\leq\hgt_M\m$ for any $\m\in\Supp_RM$, we have
\begin{eqnarray*}
\cograde_{R}(J(M),M)&=&\min\{\cograde_{R}(\m,M)~|~\m\in\Supp_RM\}\\
&\leq&\textmd{max}\{\hgt_M\m~|~\m\in\Supp_RM\}\\
&\leq&\Cdim_RM,
\end{eqnarray*}
where the third inequality follows from Lemma \ref{CosSupp}. Hence, by definition of co-Cohen Macaulay $R$-module,
we get $\hgt_M\m=\Cdim_RM$ for any $\m\in\Supp_RM$.
\end{proof}

\begin{Proposition}\label{Prop:CoCM-Catenary}
Let $R$ be a $U$ ring, $M$ a co-Cohen Macaulay $R$-module. If $\p,\q\in\Cos_RM$, with $\q\subseteq\p$. Then all saturated chains of prime ideals starting from $\q$ and ending at $\p$ have the same finite length equal to $\hgt(\p/\q)$.
\end{Proposition}

\begin{proof} Suppose that $\hgt(\p/\q)=s,~\hgt_M\q=n$. Let $\q=\p_{0}\subsetneq\p_{1}\subsetneq\cdots\subsetneq\p_{r}=\q$ be a saturated chain of prime ideals. By Theorem \ref{Thm:Co-CM} we have $\pi_{n}(\q,M)>0$, and $\pi_{n+r}(\p,M)>0$ follows from \cite[Corollary 5.7]{Li10}. Since $n+r=\hgt_M\q+r\leq\hgt_M\p$, then by Theorem \ref{Thm:Co-CM} we have $n+r=\hgt_M\p$. On the other hand, $r\leq s$, and $n+s\leq\hgt_M\p$. Hence $r=s$.
\end{proof}

Proposition \ref{Prop:CoCM-Catenary} show that the Co-Support of co-Cohen Macaulay modules over $U$ rings have catenary property.

\end{document}